%
\magnification=\magstep1 
\input amstex
\UseAMSsymbols          
\hoffset=0truecm \voffset=0truecm \hsize=15truecm
\NoBlackBoxes

   \font\rmk=cmr8  
   
  \font\ttk=cmtt8

\input pictex

\def\Hom{\operatorname{Hom}}
\def\d{\operatorname{d}}

\def\rad{\operatorname{rad}}

\def\add{\operatorname{add}}

\def\mod{\operatorname{mod}}
\def\End{\operatorname{End}}

\def\prodim{\operatorname{proj.dim}}

\def\arr#1#2{\arrow <1.5mm> [0.25,0.75] from #1 to #2}


  

\centerline{\bf Iyama's finiteness theorem via strongly quasi-hereditary algebras.}
		\bigskip
\centerline{Claus Michael Ringel}
		  \bigskip
{\narrower
{\bf Abstract.}
Let $\Lambda$ be an artin algebra and $X$ a finitely generated
$\Lambda$-module. Iyama
has shown that there exists a module $Y$ such that the endomorphism ring $\Gamma$
of $X\oplus Y$ is quasi-hereditary, with a heredity chain of length $n$,
and that the global dimension of $\Gamma$ is bounded by this $n$.
In general, one only knows that 
a quasi-hereditary algebra with a heredity chain of length $n$
must have global dimension at most $2n-2$. We want to show that Iyama's
better bound is related to the fact that the ring $\Gamma$ he constructs is
not only quasi-hereditary, but even left 
strongly quasi-hereditary:  By definition, the left
strongly quasi-hereditary algebras are the quasi-hereditary algebras
with all standard left modules of projective dimension at most~1.
\par}
	\bigskip
The aim of this note is to present a concise proof of Iyama's finiteness theorem.
For the benefit of the reader, it is essentially 
self-contained. Let us stress that all the main
arguments used are known: those of sections 1 and 2 are due to Iyama [I1,I2], 
section 3 follows the ideas of Auslander [A], whereas section 4 is based on our
joint work with Dlab [DR1, DR2, DR3]. Quasi-hereditary algebras with all standard left
modules of projective dimension at most 1 have been considered in various papers,
see for example [DR3], and it seems worthwhile to give them a name: 
we propose to call them left strongly quasi-hereditary. In our setting, 
the main advantage of 
working with left 
strongly quasi-hereditary, and not just quasi-hereditary algebras
lies in the fact that one avoids to deal with factor algebras of endomorphism rings. 
     \bigskip
{\bf 1\. Preliminaries.}
Let $\Lambda$ be an artin algebra. We denote by $\mod \Lambda$  the 
category of (finitely
generated left) $\Lambda$-modules. Morphisms will be written on the opposite side
of the scalars, thus if $f\:X \to Y$ and $g\:Y \to Z$ are $\Lambda$-homomorphisms
between $\Lambda$-modules, then the composition is denoted by $fg$.

Recall that the radical $\rad $ of $\mod \Lambda$ is defined
as follows: If $X,Y$ are $\Lambda$-modules and $f\:X \to Y$, then $f$ belongs to
$\rad(X,Y)$ provided for any indecomposable direct summand $X'$ of $X$ with inclusion map
$u\:X'\to X$ and any indecomposable direct summand $Y'$ of $Y$ with projection map $p\:Y \to
Y'$, the composition $ufp\:X' \to Y'$ is non-invertible.

Of course, for any
$\Lambda$-module $X$, the set $\rad(X,X)$ is just the radical of the endomorphism
ring of $X$, thus 
$$
 \gamma X = X\rad(X,X)
$$ 
is the radical of $X$ when considered as a right module over its endomorphism ring,
and this is a $\Lambda$-submodule of $X$. 

    \bigskip
\vfill\eject
{\bf Proposition.} {\it Let $X$ be a $\Lambda$-module. Then
\item{\rm(1)} $X$ generates $\gamma X$.
\item{\rm(2)} Any radical map $X \to X$ factors through $\gamma X$,
\item{\rm(3)} If $X$ is non-zero, then $\gamma X$ is a proper submodule of $X$.
\item{\rm(4)} If $X = \bigoplus_i X_i$ with $\Lambda$-modules $X_i$, 
then 
$\gamma X = \bigoplus_i (X_i\cap \gamma X)$, and $X_i\cap \gamma X =
   X\rad(X,X_i).$ \par
}
	\medskip
Proof: (1) Let $\phi_1, \dots,\phi_m$
be a generating set of $\rad(X,X),$ say as a $k$-module, where $k$ is the
center of $\Lambda$. Then $\gamma X = \sum_i X\phi_i$,
thus the map $\phi = (\phi_i)_i\: X^m \to \gamma X$ is surjective.

(2) is obvious. 

(3) The ring
$\Gamma = \End(X)$ is again an artin algebra and 
the radical of a non-zero $\Gamma$-module is a proper submodule
(it is enough to know
that $\Gamma$ is semi-primary).

(4) Clearly $X\rad(X,X_i) \subseteq X_i\cap \gamma X \subseteq \gamma X.$
Thus, we only have to show that for $x\in X$ and $\phi\in \rad\End(X)$,
the element $x\phi$ belongs to $\bigoplus_i X\rad(X,X_i)$.
Let $\pi_i\:X \to X_i$ be the canonical projection, so that $y = \sum_i y\pi_i$
for all $y\in X$. Then $x\phi = \sum_i x\phi\pi_i$. But with $\phi$ also
$\phi\pi_i$ belongs to $\rad$, thus $x\phi\pi_i\in X\rad(X,X_i).$
	    \medskip
{\bf Warning.} One may be tempted to say that $X$ generates $\gamma X$
by radical maps, but this is not true!
For example, let $\Lambda$ be the path algebra of the
quiver of Dynkin type $A_2$ and $X$ the minimal projective generator (i.e. the
direct sum of the two indecomposable projective modules).
Then $\gamma X$ is simple projective and the non-zero
maps $X \to \gamma X$ are not radical maps. (What is true, is the following:
$\gamma X$ is generated by $X$ using maps which have the property that when
we compose them with the inclusion map $\gamma X \subseteq X$, then they become
radical maps.)

	\bigskip
{\bf 2\. Iteration.} We consider a fixed $\Lambda$-module $X$.
We define inductively $M_1 = X$ and $M_{t+1} = \gamma M_t$, for $t \ge 1.$
According to (3), there is some $n$ such that $M_{n+1} = 0.$
The smallest such $n$ will be denoted by $\d(X)$, and we have $\d(X) \le |X|,$
where $|X|$ denotes the length of $X$.
We define $M = \bigoplus_{i=1}^n M_i$ and $M_{>t} = \bigoplus_{i>t} M_i.$
   \medskip
{\bf Warning.} 
Note that $M_2$ usually is different from $X\rad(X,X)^2$, 
a typical example is a serial module 
with composition factors $1,1,2,1,1$ (in this order) such that
the submodule of length 2 and the factor module of length 2 are isomorphic. Here,
$X\rad(X,X)^2 = 0,$ whereas $M_2$ is simple.
	      \medskip
{\bf Proposition.}
{\it Let $i \ge 1.$ Let $N$ be an indecomposable 
direct summand of $M_i$ which is 
not a direct summand of  $M_{i+1}.$
Let 
$$
 \alpha N = M_i\rad(M_i,N).
$$ 
Then $\alpha N$ is a proper submodule of $N$ and the inclusion map
$\alpha N \to N$ is a right $M_{>i}$-approximation (and of course right minimal).}
	\medskip
Proof: First, we show that $\alpha N$ is a direct summand of $M_{i+1}.$
Namely, since $N$ is
a direct summand of $M_i$, it follows that 
$\alpha N = M_i\rad(M_i,N)$ 
is a direct summand of $M_{i+1},$ using (4) for the module $M_i$. 
Since we assume that $N$ is not a direct summand of $M_i$, we see 
that $\alpha N$ has to be a proper submodule of $N$. 

In order to see that the inclusion map $u\:\alpha N \to N$ is  
a right $M_{>i}$-approximation, 
we have to show that any map $g\:M_j \to N$
with $j > i$ factors through $u$, thus that the image of $g$ is contained in $\alpha N.$
Using inductively (1), there are natural numbers $t,t'$ and surjective maps
$$
 \left(M_i \right)^{t'} @>\eta'>> \left(M_{i+1} \right)^{t} 
  @>\eta>> M_j.
$$
We claim that the 
composition $\eta' \eta g$ 
is a radical map. Otherwise, there is an indecomposable
direct summand $U$ of $\left(M_i \right)^{t'}$ such that the composition
$$
 U @>>> \left(M_i  \right)^{t'} @>\eta' >> \left( M_{i+1} \right)^{t}
 @>\eta g>> N
$$
is an isomorphism, but then $N$ is a direct summand of $M_{i+1}$, 
which is not the case. 

It follows that the image of $\eta'\eta g$ is contained in 
$M_i\rad(M_i,N) = \alpha N.$ Since $\eta'\eta$ 
is surjective, we see that 
the image of $g$ itself is contained in $\alpha N.$ 
    \bigskip

{\bf Corollary.} 
{\it Let $N$ be an indecomposable summand of $M_i$ and of $M_j$
where $i < j$. Then $N$ is a direct summand of $M_r$ for all $i \le r \le j.$}
      \medskip
Proof: Assume that $N$ is not a direct summand of $M_{i+1}$.
Since $N$ is a direct summand
of $M_j$ and $j \ge i+1,$ we can factor the identity map
$N \to N$ through the inclusion map
$\alpha N \to N$. But then $\alpha N = N$ and $N$ is a direct summand
of $M_{i+1}$, a contradiction.
   \bigskip
Given an indecomposable direct summand $N$ of $M$, there is a unique
index $i\ge 1$ such that $N$ is a direct summand of $M_i$ but not of $M_{>i}$.
We call $i$ the {\it layer} of $N$.

   \bigskip
{\bf 3. The indecomposable projective $\Gamma$-modules.}
     \medskip
We are interested in $\Gamma = \End(M)$. Recall that the indecomposable
projective $\Gamma$-modules are of the form $\Hom(M,N)$ with $N$ an
indecomposable direct summand of $M$ and we denote by $S(N)$ the 
top of the $\Gamma$-module $\Hom(M,N)$. If $M'$ is a $\Lambda$-module,
we denote by $\Hom(M,N)/\langle M'\rangle$ the factor of $\Hom(M,N)$
modulo all maps which factor through $\add M'$.
       \medskip

{\bf Proposition.} {\it Let $N$ be an indecomposable direct summand of $M$ 
with layer $i$. 
Then the minimal right $M_{>i}$-approximation $u\:\alpha N \to N$ 
yields an exact sequence
$$
 0 @>>> \Hom(M,\alpha N) @> \Hom(M,u)>> \Hom(M,N) @>>> 
  \Hom(M,N)/\langle M_{>i}\rangle @>>> 0
$$
of $\Gamma$-modules.

{\rm (a)} The $\Gamma$-module $R(N) = \Hom(M,\alpha N)$ is a direct sum of modules
of the form $\Hom(M,N'')$ with $N''$ an indecomposable direct
summand of $M$ with layer greater than $i$.

{\rm (b)} Considering the $\Gamma$-module $\Delta(N) = \Hom(M,N)/\langle  M_{>i}\rangle$,
any composition factor of $\rad \Delta(N)$
is of the form $S(N')$ where $N'$ is an indecomposable
$\Lambda$-module with layer smaller than $i$.}
		 \medskip
Proof: Since $u$ is injective, also $\Hom(u,-)$ is injective. 
Now $\alpha N$ belongs to
$M_{i+1}$, thus $\Hom(M,\alpha N)$ is mapped 
under $u$ to a set of maps $f\:M \to N$ which
factor through a module in $\add M_{>i}$. 
But since $u$ is a right $M_{>i}$-approximation, we see
that the converse also is true: 
any map $M \to N$ which factors through a module in $\add M_{>i}$ factors
through $u$. This shows that the cokernel of $\Hom(M,u)$ 
is $\Hom(M,N)/\langle  M_{>i}
\rangle$.

Of course, $R(N)$ is projective.
If we decompose $\alpha N$ as a direct sum of
indecomposable modules $N''$, then
$\Hom(M,\alpha N)$ is a direct sum of the corresponding 
projective $\Gamma$-modules $\Hom(M,N'')$
with $N''$ indecomposable and in $\add M_{>i}$.
The layer of any indecomposable module in $\add M_{>i}$ is
greater than $i$.

Now we consider $\Delta(N)$. Let $N'$ be an indecomposable direct summand of $M$
such that $S(N')$ is a composition factor of $\Delta(N)$. This means that there is
a map $f\:N' \to N$ which does not factor through $\add M_{>i}$. 
In particular, $N'$ itself
does not belong to $\add M_{>i}$.
Assume that $N'$ belongs to $\add  M_{i}$.
Also $N$ is in $\add M_i$ and according to 
(2), any radical map $M_i \to M_i$ factors through 
$M_{i+1}$. This shows that $f$ has to be invertible and therefore 
we deal with the top composition factor of $\Delta(N)$. 
It follows that the composition factors of 
$\rad\Delta(N)$
are of the form $S(N')$ with $N'$ indecomposable with layer smaller than
$i$. 
     \bigskip
{\bf 4\. Left strongly quasi-hereditary algebras.}
Let $\Gamma$ be an artin algebra. Let 
$\Cal S = \Cal S(\Gamma)$ be the
set of isomorphism classes of simple $\Gamma$-modules.
For any module $M$, let $P(M)$ be the projective cover of $M.$
 
We say that $\Gamma$ is {\it left strongly quasi-hereditary with $n$
layers} provided there is 
a function $l\:\Cal S \to \{1,2,\dots,n\}$ (the layer function)
such that for any $S \in \Cal S$, there is an exact sequence
$$
 0 \to R(S) \to P(S) \to \Delta(S) \to 0
$$
with the following two properties: 
\item{(a)}  $R(S)$ is a direct sum of projective
modules $P(S'')$ with $l(S'') > l(S)$, and
\item{(b)} if $S'$ is a composition factor of
$\rad \Delta(S),$ then $l(S') < l(S)$. 
      \medskip
Recall that $\Gamma$ is said to be {\it quasi-hereditary} with respect to
a function $l\:\Cal S \to \{1,2,\dots,n\}$ provided 
for any $S \in \Cal S$, there is an exact sequence
$$
  R(S) \to P(S) \to \Delta(S) \to 0
$$ 
with properties (a) and (b) mentioned above and the additional property
\item{(c)} For any $S\in \Cal S$, the module $P(S)$ has a $\Delta$-filtration 
(i.e. a filtration with
factors of the form $\Delta(S')$ with $S'\in \Cal S$). \par
\noindent
(But note that here the map $R(S) \to P(S)$ is not required to be injective.)
     \bigskip
{\bf Proposition.} {\it If $\Gamma$ is left 
strongly quasi-hereditary with $n$ layers and layer function $l$, then $\Gamma$
is quasi-hereditary with respect to $l$, 
and the global dimension of $\Gamma$ is at most $n$.}
    \medskip
Proof. In order to see that $\Gamma$ is quasi-hereditary with respect to $l$, 
we have to verify property (c). This we show by
decreasing induction on $l(S)$. If $l(S) = n$, then $P(S) = \Delta(S).$
Assume we know that all $P(S)$ with $l(S) > i$ have a $\Delta$-filtration.
Let $l(S) = i$. Then $R(S)$ is a direct sum of projective
modules $P(S')$ with $l(S') > l(S)$, thus it has a $\Delta$-filtration. Then also
$P(S)$ has a $\Delta$-filtration. This shows that $\Gamma$ is quasi-hereditary
with respect to $l$.

Now we have to see that the global dimension of $\Gamma$ is at most $n$. We show by induction on
$l(S)$ that $\prodim S \le l(S).$ We start with $l(S) = 1.$ In this case, $\Delta(S) = S$,
thus there is the exact sequence $0 \to R(S) \to P(S) \to S \to 0$ with
$R(S)$ projective. This shows that $\prodim S \le 1.$
For the induction step, consider some $i \ge 2$ and 
assume that $\prodim S' \le l(S')$ for all $S'$ with $l(S') < i.$ Let $S$
be simple with $l(S) = i$ and consider the exact sequence
$$
 0 \to R(S) \to \rad P(S) \to \rad \Delta(S) \to 0.
$$
All the composition factors $S'$ of $\rad\Delta(S)$ satisfy $l(S') < i$, thus
$\prodim S' < i.$ Also, $R(S)$ is projective, thus $\prodim R(S)= 0 < i$. This shows that
$\rad P(S)$ has a filtration whose factors have projective dimension less than $i$, 
and therefore 
$\prodim \rad P(S) < i.$ As a consequence, $\prodim S \le i.$ This completes the
induction. Since all the
simple modules have projective dimension at most $n$, the global dimension of
$\Gamma$ is bounded by $n$.
	 \bigskip
The bound for the global dimension cannot be improved in general:
{\it For $n\ge 2$, there are left 
strongly quasi-hereditary algebras $\Gamma$ with $n$ layers such that the global dimension of 
$\Gamma$ is equal to $n$.} As an example, take the cyclic quiver with vertices $1,2,\dots,n$, arrows
$\alpha_i\:i \to i\!-\!1$ (modulo $n$) and with relations $\alpha_{i-1}\alpha_i = 0$
for $2\le i \le n.$ 
The indecomposable projective modules $P(i)$ have the following shape:
$$ 
\hbox{\beginpicture
\setcoordinatesystem units <1cm,1cm>
\put{\beginpicture
\setcoordinatesystem units <.2cm,.25cm>
\put{$\ssize 1$} at 0 2
\put{$\ssize n$} at 0 1
\put{$\ssize n-1$} at 0 0
\endpicture} at 0 0 

\put{\beginpicture
\setcoordinatesystem units <.2cm,.25cm>
\put{$\ssize 2$} at 0 2
\put{$\ssize 1$} at 0 1
\put{} at 0 0
\endpicture} at 1 0 

\put{\beginpicture
\setcoordinatesystem units <.2cm,.25cm>
\put{$\ssize 3$} at 0 2
\put{$\ssize 2$} at 0 1
\put{} at 0 0
\endpicture} at 2 0 
\put{$\cdots$} at 3 0 

\put{\beginpicture
\setcoordinatesystem units <.2cm,.25cm>
\put{$\ssize n$} at 0 2
\put{$\ssize n-1$} at 0 1
\put{} at 0 0
\endpicture} at 4 0

\endpicture}
$$
This is a left strongly quasi-hereditary algebra using the layer
function $l(S(i)) = i$. We have $\Delta(1) = S(1)$,
whereas $\Delta(i) = P(i)$ for $2 \le i \le n.$ One easily checks that 
the projective dimension of $S(i)$ is equal to $i$, for any $i$, thus the global
dimension is $n$.

	  \bigskip
{\bf 5\. Theorem.}
{\it Let $X$ be a $\Lambda$-module. Then there is a
$\Lambda$-module $Y$ such that $\Gamma = \End(X\oplus Y)$ is left
strongly quasi-hereditary
with $\d(X)$ layers. In particular, the global dimension of $\Gamma$ is at most $\d(X)$.}
     \medskip
In addition, we record that $\d(X) \le |X|$. Also, the  
construction of $Y$ shows that we can assume that
{\it any indecomposable direct summand of the module $Y$ is
a submodule of an indecomposable direct summand of $X$.} 
  \medskip
Proof: As above, let $M_1 = X$ and $M_{t+1} = \gamma M_t$ for $t\ge 0.$ 
Let $Y = M_{>1}$ and $M = X\oplus Y = \bigoplus_{i=1}^n M_i$ with $n = \d(X)$.
According to Proposition 3, $\Gamma$ is a left strongly quasi-hereditary algebra
with $n$ layers, thus we can apply Proposition 4.

The additional information comes from (3) and (4).
    \bigskip
Several applications should be mentioned.
	\medskip
{\bf (1)} First, there is 
the representation dimension of an artin algebra $\Lambda$ as introduced by Auslander
in the Queen Mary Notes.
By definition, this is the smallest number which occurs as the
global dimension of the endomorphism ring of a $\Lambda$-module which is 
both a generator and a cogenerator. If one takes $X = \Lambda \oplus D\Lambda$, where
$D = \Hom_k(-,k)$ is the $k$-duality functor, then 
the theorem provides a $\Lambda$-module $Y$ such that
the global dimension of $\End(X\oplus Y)$ is bounded by $|X| = 2|\Lambda|.$
Since $X\oplus Y$ is a generator and a cogenerator, this yields a bound
for the representation dimension of $\Lambda.$ In particular, {\it the representation
dimension is always finite,} this is Iyama's  finiteness theorem.                           
	  \medskip
{\bf (2)} In this way, we obtain for artin algebras
also a strenghtening of the main result of [DR2]: {\it If $\Lambda$
is an artin algebra, then there is an artin algebra $\Gamma$ and an idempotent $e\in 
\Gamma$ with $e\Gamma e = \Lambda$ such that $\Gamma$ is  left strongly quasi-hereditary,} 
and not only quasi-hereditary.
Here, one may start with $X = \Lambda$
(or also with any module which has $\Lambda$ as a direct summand, 
as the one in the previous paragraph).
   \bigskip
{\bf (3)} Finally, we see that Auslander algebras are 
left strongly quasi-hereditary: {\it If $\Lambda$ is a
representation-finite artin algebra and $M$ is the direct sum of the
indecomposable $\Lambda$-modules, one from each isomorphism class, then
$\Gamma = \End(M)$ is left strongly quasi-hereditary.} Namely, the Theorem
asserts that there is a $\Lambda$-module $Y$, such that $\End(M\oplus Y)$
is left strongly quasi-hereditary. But $\Gamma$ and $\End(M\oplus Y)$ are 
Morita equivalent, thus also $\Gamma$ is left strongly quasi-hereditary.
       \bigskip
{\bf Appendix}
    \medskip 
For the convenience of the reader, the appendix collects from
the literature some further information on left strongly quasi-hereditary
algebras. Also, we will add some examples which may be useful.
    \bigskip 
{\bf A1. Characterizations of left strongly
quasi-hereditary algebras.}
    \medskip 
Let us assume that $\Gamma$ is quasi-hereditary with respect to
some layer function $l$. For any simple module $S$, let $\Delta(S)$ be the
corresponding standard module, $\nabla(S)$ the costandard module. Let
$T$ be the characteristic tilting module.  Given
a set $\Cal X$ of modules, we denote by $\Cal F(\Cal X)$ the class of
modules which have a filtration with all the factors in $\Cal X$.
Finally, recall that a module is
said to be {\it divisible} provided it is generated by an injective module.
    \medskip 
\vfill\eject
{\bf Proposition.} {\it For the quasi-hereditary algebra
$\Gamma$ the following conditions are equivalent: \item{\rm(1)} Any
$\Delta$-module has projective dimension at most $1$.
 \item{\rm(2)} Any module in $\Cal F(\Delta)$ has projective dimension at
most $1$.
  \item{\rm(3)} $T$ has projective dimension at most $1$. \item{\rm(4)} The
modules in $\Cal F(\nabla)$ are the modules generated by $T$.
 \item{\rm(5)} Any module generated by $T$ belongs to $\Cal F(\nabla)$.
  \item{\rm(6)} $\Cal F(\nabla)$ is closed under factor modules.
  \item{\rm(7)} Any divisible module belongs to $\Cal F(\nabla)$.
\item{\rm(8)} For any module $M$, there is an exact sequence
           $0\to M \to D_0\to D_1 \to 0$ where $D_0, D_1$ are modules in
$\Cal F(\nabla)$.
 \item{\rm(9)} There is an exact sequence
           $0\to \Gamma \to D_0\to D_1 \to 0$ where $D_0, D_1$ are modules
in $\Cal F(\nabla)$.\par}
    \medskip 

Before we outline the proof, let us stress the following:
Condition (3) states that $T$ is what sometimes is called a {\it classical}
tilting module, namely a tilting module of projective dimension at most $1$. 
    \medskip 
Proof. For
the equivalence of (1), (2), (6) and (7) we may refer to [DR3], Lemma 4.1
(section 5 of that paper contains also the assertion that (1) implies (4)).
Of course, (2) implies (3), and classical tilting theory asserts that (3)
implies (4). Trivially, (4) implies (5) and (6), also (6) implies (7), since
the injective modules belong to $\Cal F(\nabla)$. In order to see that (7)
implies (8), one just takes for $D_0$ the injective envelope of $M$. Again
(8) implies (9) is trivial. The equivalence of (3) and (9) is part of
tilting theory. It remains to see that (5) implies (4), but it is easy to
see that any module in $\Cal F(\nabla)$ is generated by $T$.
    \bigskip 
{\bf A2\. The missing left-right symmetry.} 
\medskip 
An artin algebra $\Gamma$ is said to be {\it
right strongly quasi-hereditary} provided the opposite algebra
$\Gamma^{\text{op}}$ is left strongly quasi-hereditary.
      \medskip 
{\bf (1)} {\it A left strongly quasi-hereditary algebra need not be
right strongly quasi-hereditary.}
      \medskip
As an example, consider the algebra $\Gamma$ with quiver
$$ 
\hbox{\beginpicture
\setcoordinatesystem units <1cm,1cm>
\put{$2$} at 0 0 
\put{$1$} at 1 0 
\put{$3$} at 2 0 
\arr{0.8 0.1}{0.2 0.1} 
\arr{0.2 -.1}{0.8 -.1} 
\arr{1.2 0.1}{1.8 0.1} 
\arr{1.8 -.1}{1.2 -.1} 

\put{$\alpha$}  at 0.5  0.4
\put{$\alpha'$} at 0.5  -.4
\put{$\beta$}   at 1.5  0.4 
\put{$\beta'$}  at 1.5  -.4
\endpicture}
$$
and with relations $\alpha\alpha', \beta\alpha', \beta\beta', \alpha'\alpha\beta'$. 
The indecomposable projective modules $P(i)$ have the following shape:
$$ 
\hbox{\beginpicture
\setcoordinatesystem units <1cm,1cm>
\put{\beginpicture
\setcoordinatesystem units <.2cm,.25cm>
\put{$\ssize 1$} at 1 3
\put{$\ssize 2$} at 0 2
\put{$\ssize 1$} at 0 1
\put{$\ssize 3$} at 2 2
\put{$\ssize 1$} at 2 1
\put{$\ssize 2$} at 2 0
\endpicture} at 0 0 
\put{\beginpicture
\setcoordinatesystem units <.2cm,.25cm>
\put{$\ssize 2$} at 1 3
\put{$\ssize 1$} at 1 2
\endpicture} at 1,3 0.22 
\put{\beginpicture
\setcoordinatesystem units <.2cm,.25cm>
\put{$\ssize 3$} at 1 3
\put{$\ssize 1$} at 1 2
\put{$\ssize 2$} at 1 1
\endpicture} at 2.5 0.1
\endpicture}
$$
It is obvious that the numbering of the simple modules provides a layer function
so that  $\Delta(2)$ and $\Delta(3)$ are projective, whereas $\Delta(1)$ is simple
with an exact sequence
$$
 0 \to P(2)\oplus P(3) \to P(1) \to \Delta(1) \to 0.
$$
Instead of looking at modules over the opposite algebra, we can consider their
$k$-duals. If $\Gamma^{\text{op}}$ would be right strongly quasi-hereditary,
we would obtain an exact sequence of $\Gamma$-modules of the form
$$
 0 \to \nabla(1) \to I(1) \to Q(1) \to 0,
$$
where $I(1)$ is the injective envelope of $1$, where $\nabla(1)$ has only one
composition factor of the form $1$ and where $Q(1)$ is injective.
The indecomposable injective modules  have the shape
$$ 
\hbox{\beginpicture
\setcoordinatesystem units <1cm,1cm>
\put{\beginpicture
\setcoordinatesystem units <.2cm,.25cm>
\put{$\ssize 1$} at 0 2
\put{$\ssize 2$} at 0 1
\put{$\ssize 1$} at 2 2
\put{$\ssize 3$} at 2 1
\put{$\ssize 1$} at 1 0
\endpicture} at 0 0 
\put{\beginpicture
\setcoordinatesystem units <.2cm,.25cm>
\put{$\ssize 1$} at 1 3
\put{$\ssize 3$} at 1 2
\put{$\ssize 1$} at 1 1
\put{$\ssize 2$} at 1 0
\endpicture} at 1,3 0.1 
\put{\beginpicture
\setcoordinatesystem units <.2cm,.25cm>
\put{$\ssize 1$} at 1 1
\put{$\ssize 3$} at 1 0
\endpicture} at 2.5 -0.1
\endpicture}
$$
Since $I(1)$ contains three composition factors of the form $1$, 
we see that $Q(1) \neq 0$,
thus $\nabla(1)$ has injective dimension equal to $1$.
But the only submodule of $I(1)$ with injective dimension equal to $1$ is of length 3
with two composition factors $1$ (and one composition factor $2$).
This shows that $\Gamma$ cannot be  right strongly quasi-hereditary. 
     \bigskip

{\bf (2)} Let $\Gamma$ be left strongly quasi-hereditary with layer function
$l\:\Cal S(\Gamma) \to \{1,2,\dots,n\}$. Again, let $D = \Hom_k(-,k)$. If
$S$ is a simple $\Gamma$-module, then we define $l(DS) = l(S)$, thus we
consider $l$ also as a function $l\:\Cal S(\Gamma^{\text{op}}) \to
\{1,2,\dots,n\}$. We have shown above that $\Gamma$ is quasi-hereditary with
respect to $l$, and it is well-known that then also $\Gamma^{\text{op}}$ is
quasi-hereditary with respect to $l$. In general, {\it $\Gamma^{\text{op}}$
may not be left strongly quasi-hereditary with respect to $l$, even if it is
left strongly quasi-hereditary with respect to some other layer function.}
     \medskip 
As a typical example, consider an algebra $\Gamma$ such that
the quiver of $\Gamma$ has no oriented cycles. Then there is a layer
function $l$ such that the $\Delta$-modules are projective, and then the
standard modules for the opposite algebra are the simple modules. In this
case, $\Gamma$ is left strongly quasi-hereditary with respect to $l$, but it
is right strongly quasi-hereditary with respect to this $l$ only in case
$\Gamma$ is hereditary. But of course, always $\Gamma^{\text{op}}$ will be
left strongly quasi-hereditary, however we have to use a different layer
function. 

Also, the example exhibited in section 4 is of this kind: The
algebra $\Gamma^{\text{op}}$ is not left strongly quasi-hereditary with
respect the the ordering $\{1,2,\dots,n\}$, but it is left strongly
quasi-hereditary with respect to the ordering $\{n-1,n,1,\dots,n-2\}$.
    \bigskip 

In fact, there is the following general result due to
Erdmann-Parker ([EP],2.1):
    \medskip 
{\bf Proposition.} {\it If $\Gamma$ is both left strongly
quasi-hereditary and right strongly quasi-hereditary with respect to the same function
$l$, then the global dimension of $\Gamma$ is at most $2$.}
    \medskip 
Proof: The implication (1) $\implies$ (7) of Propositon A1
for $\Gamma^{\text{op}}$
shows that if $\Gamma$ is right strongly quasi-hereditary, then all
submodules of projective modules belong to $\Cal F(\Delta)$. If $\Gamma$ is
left strongly quasi-hereditary, then the modules in $\Cal F(\Delta)$ have
projective dimension at most $1$. But if all submodules of projective
modules have projective dimension at most $1$, then the global dimension of
$\Gamma$ is at most $2$.
    \bigskip 
{\bf (3)} 
If $\Gamma$ is quasi-hereditary with characteristic tilting
module $T$, then the endomorphism ring $\Gamma'$ of $T$ is called the R-dual
(Ringel-dual) of $\Gamma$. It is again quasi-hereditary with respect to a
suitable layer function (so that the characteristic tilting module $T'$ for
$\Gamma'$ is given by $\Hom_\Gamma(T,Q)$, where $Q$ is a minimal injective
cogenerator for the category of $\Gamma$-modules). Again let us mention an
observation of Erdmann-Parker ([EP], section 3):
    \medskip 
{\bf Proposition.} {\it The R-dual of a left strongly
quasi-hereditary algebra is right strongly quasi-hereditary.}
    \medskip 
Proof: Tilting theory asserts: if $T$ is a
tilting module of projective dimension 1, then the injective dimension of
$T' = \Hom_\Gamma(T,Q)$ is at most $1$.
    \bigskip
{\bf A3. Historical remarks.}
    \medskip 
Left strongly quasi-hereditary algebras have been considered in
various papers, only the name is new. As we have mentioned, several
characterizations of these algebras have been given already in 1992 in our
joint survey [DR3] with Dlab.
    \medskip 
It is obvious that any hereditary artin algebra is left
strongly quasi-hereditary with respect to any total ordering of the simple
modules (thus quasi-hereditary with respect to any total ordering [DR1]).
    \medskip 
Other important examples of left strongly quasi-hereditary algebras
are the Auslander algebras.
Several papers by Br\"ustle, Hille and R\"ohrle, but also others are devoted to such
examples.
    \medskip 
Under suitable directedness assumptions, bimodule problems can
be described using left strongly quasi-hereditary algebras, see for example
Hille and Vossieck [HV].
    \medskip 
In the context of preprojective algebras, Geiss, Leclerc and
Schr\"oer [GLS] have shown that endomorphism rings of suitable rigid modules
are left strongly quasi-hereditary, and this has been generalized by Iyama
and Reiten [IR] and by Schr\"oer himself [S].
    \medskip 
On the other hand, in the classical realm of the quasi-hereditary arising for
semisimple Lie algebras and algebraic groups, one cannot expect that the
quasi-hereditary algebras occurring there are left strongly quasi-hereditary.
The reason is quite simple: Usually, these quasi-hereditary algebras are R-self-dual, 
and have quite large global dimension. However, R-self-dual
algebras which are left
strongly quasi-hereditary are also right strongly quasi-hereditary,
and thus they have global dimension at most 2 (see the results of Erdmann-Parker
mentioned in A2). 
	  \medskip
But note that in this setting, the equivalence of the
conditions (1) and (8) mentioned in A1 was already formulated by Friedlander
and Parshall ([FP, Proposition 3.4]), before the concept of a quasi-hereditary algebra
was introduced. Following [FP] one may say that the
$\nabla$-filtration dimension of a module $X$ is at most $d$ provided there
exists an exact sequence 
$$
 0 \to X \to D_0 \to D_1 \to \cdots \to D_d \to 0 
$$ 
with $D_0,\dots, D_d\in\Cal F(\nabla).$ 
Using this terminology,
 the equivalence of (1) and (8) may be reformulated as the following
assertion: {\it A quasi-hereditary algebras $\Gamma$ is left strongly
quasi-hereditary if and only if the global $\nabla$-filtration dimension of
$\Gamma$ is at most $1$.}
    \bigskip\bigskip 
{\bf Acknowledgment.} This note is based on lectures given at Shanghai
and Bielefeld in spring 2008. The author is grateful for many helpful comments by
the audience.
In addition he wants to thank the referee for a careful reading of the paper.

   \bigskip\bigskip
{\bf References.}
\frenchspacing
	\medskip
\item{[A]}
Auslander, M.: The representation dimension of artin algebras. Queen Mary
College Mathematics Notes (1971)

\item{[DR1]}
Dlab, V., Ringel C. M.: Auslander algebras as quasi-hereditary algebras. 
J. London Math. Soc. 39 (1989), 457-466. 

\item{[DR2]}
Dlab, V., Ringel C. M.: Every semiprimary ring is the endomorphism ring of a projective 
module over a quasi-hereditary ring. Proc. Amer. Math. Soc. 107 (1989), 1-5.

\item{[DR3]}
Dlab, V., Ringel C. M.: The module theoretical approach to
quasi-hereditary algebras. In: Representations of Algebras and Related Topics.
London Math. Soc. Lecture Note Series 168 (1992), 200-224.

\item{[EP]} Erdmann, K,
Parker, A.: On the global and $\nabla$-filtration dimensions of
quasi-hereditary algebras. J Pure Appl. Algebra 194 (2004), 95-111.

\item{[FP]} Friedlander, E.M., Parshall, B.: Cohomology of Lie algebras and
algebraic groups. Amer.J.Math. 108 (1986), 235-253. 

\item{[GLS]} Geiss, C., Leclerc, B, Schr\"oer, J.: 
 Cluster algebra structures and semicanonical bases for unipotent groups.
Preprint. arXiv:math/0703039 

\item{[HV]} Hille, L., Vossieck, D.:
 The quasi-hereditary algebra associated to the radical bimodule over a
 hereditary algebra, Colloq. Math. 98 (2003), 201.211.

\item{[I1]}
Iyama, O.: Finiteness of representation dimension. Proc\. Amer\. Math\. Soc.
131 (2003), 1011-1014.

\item{[I2]}
Iyama, O.: Rejective subcategories of artin algebras and orders.
\newline arXiv:math/0311281. (Theorem 2.2.2 and Theorem 2.5.1).

\item{[IR]} Iyama, O., Reiten, I.:
2-Auslander algebras associated with reduced words in Coxeter group. In preparation.

\item{[S]} Schr\"oer, J.: Flag varieties and quiver Grassmannians. Lecture at the
Conference: Homological and geometric methods in algebra. Trondheim 2009.

                                                  \bigskip\bigskip
{\rmk Fakult\"at f\"ur Mathematik, Universit\"at Bielefeld,
POBox 100\,131, \ 
D-33\,501 Bielefeld}

{\rmk E-mail address:} {\ttk ringel\@math.uni-bielefeld.de}

\bye